# Mathematical Models for Optimization of Grid-Integrated Energy Storage Systems

Chiebuka Eyisi[*], *Student Member IEEE*, Ameena Saad Al-Sumaiti[β], *Member IEEE*, Konstantin Turitsyn[ψ], *Member IEEE*, and Qifeng Li[*], *Member IEEE*,

[*] Dept. of Elect. & Comp. Engr., University of Central Florida, Orlando, FL, USA (e-mail: cvpeyisi@knights.ucf.edu; qifeng.li@ucf.edu)
[β] Dept. of Elect. & Comp. Engr., Khalifa University, Abu Dhabi, UAE (e-mail: ameena.alsumaiti@ku.ac.ae)
[ψ] Dept. of Mech. Engr., Massachusetts Institute of Technology, Cambridge, MA, USA (e-mail: turitsyn@mit.edu)

*Abstract*—Energy storage has been proven to yield positive effects on planning, operation and control of electric grids. It has become a crucial task to properly model the energy storage systems (ESS) under the framework of grid optimization on transmission and distribution networks including microgrids. This paper presents a review on mathematical models and test cases of ESSs used for grid optimization studies, where the network constraints of power systems are included. The existing ESS models are mainly classified into two categories – linear and nonlinear models. The two main categories are further divided into several sub-categories respectively; such as mixed-integer linear and convex nonlinear sub-categories. Based on the review and discussions, this paper aims at providing suggestions for choosing proper ESS models for specific grid optimization studies considering the chosen power network model.

## I. Introduction

Energy storage has attracted substantial attentions worldwide and will continue to attract more in the future [1]. A major objective for many electric power utilities is to gain more flexibility in managing the supply and demand of electric power in order to enhance the power system's operational stability and efficiency. The technology of Energy Storage Systems (ESS) refers to the process of converting energy from one form (mainly electrical energy) to a storable form (during off-peak hours) and reserving it in various mediums; then the stored energy can be converted back into electrical energy when needed (during peak hours) [2]. Recent developments in ESS technologies have helped mitigate the imbalance between power supply and demand to a degree and enabled higher penetrations of Renewable Energy Sources (RES). In the last few years, research interests have been fueled into the capability of grid integrated ESS to provide support for the grid and smoothen the output of RES [3].

In some existing literature, the authors have reviewed and summarized ESS technologies, their research and development, grid demonstrations and industrial applications [1], [2], [3]. However, the mathematical models of ESSs for the purposes of grid-integrated optimization studies were not well discussed in these review papers. This paper is aimed at providing the reader with detailed information as to how to model ESSs in power system optimizations. The focus are implementations in transmission and distribution networks including microgrids. The coordination between the ESS models and power network flow models is also discussed.

The rest of the paper is organized as follows: Section II presents a brief summary on ESS technologies that are suitable for power grid applications. Section III lists and compares various mathematical models for grid-integrated ESS and recommends models suitable for an ESS technology. In Section IV, information regarding test systems as seen in literature is provided, before making concluding remarks in Section V.

## II. ESS Technologies For Grid Applications

### A. Types of ESS

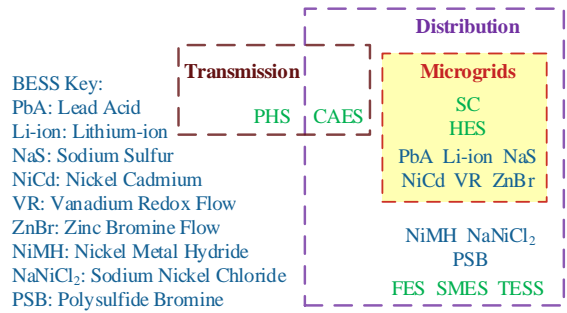

Figure 1. ESS Classification

Storage elements for electrical grids include, but are not limited to: Pumped Hydro Storage (PHS) [4], Compressed Air Energy Storage (CAES) [5], Flywheel Energy Storage (FES) [6], Battery Energy Storage Systems (BESS) [7], Super Capacitors (SC) [8], Thermal Energy Storage Systems (TESS) [9], Superconducting Magnetic Energy Storage (SMES) [10] and Hydrogen Energy Storage (HES) [11]. These ESS technologies can be categorized in terms of their suitable storage durations, response times, functions or the stored form of energy. Here we categorize these ESS technologies as shown in Figure 1 above with the aid of [1] – [11], based on network deployment levels.

C. Eyisi and Q. Li's work was supported by U.S. National Science Foundation under award 1808988.

A. S. Al-Sumaiti and K. Turitsyn's work was supported by the MI-MIT Cooperative Program under grant MM2017-000002.

## B. Functions of ESS in Power Grids

Some of the benefits of ESS in electric power grids include energy management, load leveling, voltage support, alleviating the intermittence of RES, frequency regulation, potential for financial profits, providing backup power during emergencies, grid stabilization and resiliency [2], [3]. We identify potential for financial profits, reliability and stability benefits as the leading functions of ESS in power grids. In today's electricity markets, the profit-maximizing size of an ESS is primarily determined by its technological characteristics (self-discharge rate and round-trip charge/discharge efficiency). This may provide revenue opportunities to RES owners for incorporating ESSs engaging in price arbitrage [12].

Power system reliability involves the considerations of the states of the system and whether they are adequate, secure or insecure. Integrating ESS within a power system ameliorates reliability challenges; by helping meet the demand when RES or conventional generation sources are insufficient, or prolonging the lifetime of these sources through more frequent scheduled maintenance, or reconfiguring a part of the network capable of operating in islanded mode with ESS installations, and helping re-dispatch generation sources in the network during power outages lasting more than a few seconds. Reliability indices like Loss of Load Expectation (LOLE) and Loss of Load Probability (LOLP) are typical indices used in the reliability evaluation of power systems with ESS [13].

In maintaining power system stability, load balance can be managed through conventional generators in conjunction with RES and ESS. ESS being able to inject and absorb real power helps reduce power oscillations. ESS can also help maintain power system stability and necessary voltage levels through reactive power support. The combination of Flexible AC Transmission Systems (FACTS) and ESS can help to improve system stability and power transfer capability of the network. Investigations related to ESS for power system stability can be found in [14].

The aforementioned ESS functions are realizable when discussing the concept of a microgrid. A microgrid is a group of distributed RESs and interconnected loads capable of acting as a single controllable entity within defined electrical boundaries with respect to the external grid. Its characteristic attribute is its capability to disconnect (operating in the islanded mode) and reconnect to the external grid. With ESS being able to stabilize these systems and enable integration of RES thereby improving power system reliability, microgrid concepts are realizable from the viewpoints of planning, operation and control [15].

## III. MATHEMATICAL MODELS

During the planning and operation stages of an ESS for distribution and transmission networks, questions of size, site and costs needs to be addressed. In the power systems industry, optimization models help in the decision-making process of scheduling and dispatching power generation resources. Formulations such as Optimal Power Flow (OPF) [16] and Unit Commitment (UC) [17] represents the majority of these descriptive models. In general, the power system optimization considering ESSs has the following form:

$$\min \quad objective\ function \quad (1a)$$
$$s.t. \quad network\ constraints \quad (1b)$$
$$\quad energy\ storage\ constraints \quad (1c)$$

The objective function in the above model can be any but not limited to the following: active power cost minimization, active and reactive power loss minimization, emissions minimization, voltage deviation minimization, power transfer maximization and several others [18]. These objective functions are generally convex. The non-convexity and computational burden come from the type and size of the constraints in (1b) and (1c). The network constraints in (1b) have been discussed extensively in [18], [19]. In this paper, the mathematical models being referred to focus on how ESSs are modeled typically as constraints as shown in (1c).

In existing literature, some mathematical models proposed by various authors can be classified into linear and nonlinear models. Regardless of the model implemented; for planning and operational purposes, the intent is to address questions related to the size, site and associated costs of the ESS to be incorporated in the network while meeting other system reliabilities. The literature on modeling ESS for grid integration is seen to cover quite a number of indices and can be tailored to suit various planning and operational objectives for different power system optimizations. A storage device can be characterized by its rated power ($P_{ESS}$) in kW/MW, rated energy ($E_{ESS}$) in kWh/MWh and also charging and/or discharging cycle efficiencies ($\eta_{ch}$ and/or $\eta_{disch}$). As previously stated, these models are part of realizations where information about the network's topology, captured by an admittance matrix is also included and not ignored. Hence power flow constraints related to the network's topology are included in the model to be optimized. We begin by showing what these ESS models might constitute in implementations.

### A. Linear Models
#### 1) A Mixed Integer Linear Model:

In the following model, charge and discharge efficiencies are considered while binary variables are introduced to formulate the charge/discharge state of the ESS.

$$P_{ESS,n,t}^{net} = \eta_{ch} P_{ESS,n,t}^{ch} - (1/\eta_{disch}) P_{ESS,n,t}^{disch} \quad (2a)$$

$$E_{ESS,n}^{min} \leq E_{ESS,n,0} + \Delta t \sum_{k=1}^{t} P_{ESS,n,k}^{net} \leq E_{ESS,n}^{max} \quad (2b)$$

$$0 \leq P_{ESS,n,t}^{ch} \leq \alpha_{n,t}^{ch} P_{ESS,n,t}^{max} \quad (2c)$$

$$0 \leq P_{ESS,n,t}^{disch} \leq \alpha_{n,t}^{disch} P_{ESS,n,t}^{max} \quad (2d)$$

$$\alpha_{n,t}^{ch} + \alpha_{n,t}^{disch} \leq 1 \quad (2e)$$

$$\alpha_{n,t}^{ch}, \alpha_{n,t}^{disch} \in \{0,1\},\ n \in N,\ t = 0,1,\dots,T-1 \quad (2f)$$

where (2a) describes the net injected/absorbed power to the ESS, $P^{net}$ using the charge ($P^{ch}$) and discharge ($P^{disch}$) variables at each ESS connected bus $n$ and at time $t$. It is also dependent on the charge and discharge cycles' efficiencies. The charging efficiency should be less than 1 because a percentage of energy fed to the ESS will be stored while the remaining will be lost in the form of losses. Conversely, to withdraw energy from the

ESS, more energy is required to cover the discharging losses and hence this is the reason for $1/\eta_{disch}$ in (2a). Equation (2b) ensures the current available amount of stored energy in the ESS is within its technological minimum and maximum, where $E_{ESS,n,0}$ is the initial energy stored at the beginning of the interval at bus $n$. The parameter $\Delta t$ is the time step resolution in hours (for example, $\Delta t$ is 0.25 hours if 15-minutes data is used). Equations (2a) and (2b) are called State of Charge (SOC) constraints. Equations (2c) and (2d) ensure the charge and discharge powers do not exceed the ratings of the ESS. The charge and discharge indicator binary variables $\alpha^{ch}$ and $\alpha^{disch}$ respectively in (2c), (2d), (2e) restricts the ESS from charging and discharging at the same time $t$. The index $T$ is the number of periods under study and $N$ is the set of all buses with an ESS connected. These constraints represent the majority of constraints seen in literature. Components of the above model were illustrated in [20] – [25].

*2) A Continuous Linear Model:*

To eliminate the binary variables, some researchers used the following simplified lossless model.

$$-P_{ESS,n}^{max} \leq P_{ESS,n,t} \leq P_{ESS,n}^{max} \quad (3a)$$

$$E_{ESS,n}^{min} \leq E_{ESS,n,0} + \Delta t \sum_{k=1}^{t} P_{ESS,n,k} \leq E_{ESS,n}^{max} \quad (3b)$$

$$n \in N, \ t = 0, 1, \dots, T-1 \quad (3c)$$

where the power ($P_{ESS}$) flowing into the ESS located at bus $n$ at each time $t$ can be either positive (discharging) or negative (charging) and is bound by the constraint in (3a). The constraint in (3b) bounds the maximum capacity of each ESS till any time period less than the planning horizon by taking into consideration the initial energy stored at the beginning of the interval, $E_{ESS,n,0}$. $N$ is the set of all buses with an ESS connected and $T$ is the planning horizon. Components of the above model were illustrated in [26], [27].

## B. Nonlinear Models

*1) Nonlinear Model Adapted from Linear Models by Replacing Integer Indicator Variables with a Product of the Control Variables:*

These models are adapted from any of the linear models in which a product of any of the control variables related to ESS overrides the application of integer variables used to indicate and avoid simultaneous charging and discharging within the same time period.

$$P_{ESS,n,t}^{ch} \cdot P_{ESS,n,t}^{disch} = 0 \quad (4a)$$

$$0 \leq P_{ESS,n,t}^{ch} \leq P_{ESS,n,t}^{max} \quad (4b)$$

$$0 \leq P_{ESS,n,t}^{disch} \leq P_{ESS,n,t}^{max} \quad (4c)$$

where (3c) and the SOC constraints given in (2a) and (2b) are also applicable here. The equation (4a) guarantees that at most one of the charge/discharge variables is non-zero. Equations (4b) and (4c) are similar to (2c) and (2d) respectively. The above model is illustrated in [28], [29]. An interesting method of removing the bilinear constraint (4a) by adding a disturbance term to the objective function was provided in [29].

*2) Nonlinear Model for BESS:*

This nonlinear model was developed based on an equivalent electric circuit of a BESS (see Fig. 4 in [30]).

$$P_{ESS,n,t}^{loss} V_{n,t} = r_n^{eq}(P_{ESS,n,t})^2 + r_n^{CVT}(Q_{ESS,n,t})^2 \quad (5a)$$

$$P_{ESS,n,t}^{net} = P_{ESS,n,t} + P_{ESS,n,t}^{loss} \quad (5b)$$

$$(P_{ESS,n,t})^2 + (Q_{ESS,n,t})^2 \leq (S_{CVT,n}^{max})^2 \quad (5c)$$

where (3c) and the SOC constraint in (2b) are applicable here with $P^{net}$ now incorporating the ohmic power losses ($P^{loss}$) of the BESS. The equation (5a) relates the ohmic power losses of the BESS in terms of its active and reactive power outputs, where the equivalent resistance, $r^{eq} = r^{BESS} + r^{CVT}$; is the sum of the BESS and associated converter resistances. The power ($P_{ESS}$) flowing into the BESS located at bus $n$ and at each time $t$ can be either positive (discharging) or negative (charging). The variable, ($V_{n,t}$) in (5a) is the square of the voltage magnitude at bus $n$. The active and reactive power limits of the ESS is set to not exceed the MVA ratings of the converter in (5c). The above model is illustrated in [30] and similar approximations can be found in [31], [32]. The authors in [33] also presented approximations of the above model by proposing a combined active-reactive optimal power flow in distribution networks using charge/discharge dynamics that incorporated similar SOC deductions as in (2a) and (2b).

*3) Convex Nonlinear Model for BESS:*

This model, similar to the above nonlinear model in (5) for BESS is shown in [34]; where the authors presented a novel convex relaxation of the non-convex equation (5a) to their distributed energy storage optimal scheduling problem.

$$P_{ESS,n,t}^{loss} V_{n,t} \geq r_n^{eq}(P_{ESS,n,t})^2 + r_n^{CVT}(Q_{ESS,n,t})^2 \quad (6a)$$

$$r_n^{eq}(S_{CVT,n}^{max})^2 \geq P_{ESS,n,t}^{loss} V_n^{min} + r_n^{BESS}(Q_{ESS,n,t})^2 \quad (6b)$$

$$(S_{CVT,n}^{max})^2 (V_n^{min} + V_n^{max})$$
$$\geq (S_{CVT,n}^{max})^2 V_{n,t} + P_{ESS,n,t}^{loss}(V_n^{min} \cdot V_n^{max}) \quad (6c)$$

where (3c), the SOC constraint in (2b), (5b) and (5c) are also applicable here. The equations (6a – 6c) is the convex hull relaxation of (5a) within system bounds ($V_n^{min} \leq V_{n,t} \leq V_n^{max}$ and (5c)). Please refer to [34] for more information on this model.

## C. Discussion

This subsection discusses suggestions on selecting a proper ESS model for constraint (1c) when the network constraints in (1b) is defined.

*1) Transmission Networks:*

*a) If nonlinear AC network model is used in (1b):*

We suggest models (3) or (4) be chosen so as to avoid yielding a Mixed Integer Nonlinear Problem (MINLP) which is difficult to solve. If model (3) or (4) is chosen, the resulting power system optimization problem in (1) is a Nonlinear Problem (NLP), which is relatively easier to solve than an MINLP. Model (4) is more accurate than model (3); however, introduces more computational burden.

*b) If linearized AC network model is used in (1b):*

We suggest models (2) or (3) be chosen; again to avoid yielding an MINLP. The linearized AC network model in (1b) incorporates certain reasonable assumptions peculiar to transmission networks. If model (3) is chosen, the power system optimization problem in (1) is a Linear Problem (LP) which is easier to solve, whereas if model (2) is chosen, (1) becomes a Mixed Integer Linear Problem (MILP) which although is harder than an LP, is more accurate.

*2) Distribution Networks and Microgrids:*

For distribution networks, careful attention has to be paid to them due to their low reactance/resistance (X/R) ratios and their radial or weakly meshed structures. Linearizing the nonlinear equations in (1b) for distribution networks is not ideal, hence an AC network model should be used. Any of the non-convex models (2), (3) or (4) is recommended to be paired with (1b) for grid optimization in distribution networks. If model (3) or (4) is chosen, the resulting power system optimization problem in (1) is an NLP, whereas if model (2) is used, (1) becomes an MINLP which is harder than an NLP, but more accurate.

For BESSs, we recommend model (5). The more nonlinear equations present in a power system optimization formulation, the more difficult and complex it is to solve. Convexification is now a practical way to reduce computational complexity for solving a non-convex problem. Model (6) provides an opportunity to researchers interested in the optimization of grid-integrated BESSs using convex techniques.

TABLE I.  CLASSIFYING MATHEMATICAL MODELS

| Model | Incorporated By |
|---|---|
| (2) | [12], [20]$^d$, [21], [22]$^d$, [23], [24], [25], [37], [38], [39]$^d$, [40]$^d$, [41], [42], [43], [44]$^d$, [45] |
| (3) | [26], [27] |
| (4) | [28]$^d$, [29]$^d$ |
| (5) | [30]$^d$, [31]$^d$, [32]$^d$, [33] |
| (6) | [34]$^d$ |

d – applied to distribution network

In addition, some authors may have components common to any of the five models presented. This is common practice, as researchers, planning engineers and power system operators might need to capture several aspects from different models to suit their problem. A literature on grid-integration of ESS technologies can be seen in [35], [36]. Table I attempts to classify some of the applicable grid-integrated linear and nonlinear ESS models used in power system optimization.

TABLE II.  MODELS FOR ESS

| ESS | Related References |
|---|---|
| PHS | [25], [42], [43] |
| CAES | [38], [43] |
| FES | [45] |
| BESS | [21], [23], [24],[30], [31], [32],[33], [34], [37], [39] |
| SMES | Generic {[12], [20], [26], [27], [28], [29], [41]} |
| SC | Generic {[12], [20], [26], [27], [28], [29], [41]} |
| TESS | [22], [45] |
| HES | [40] |

Table II is intended to provide the reader with suitable models incorporating other technological aspects and characteristics not covered in this paper for all the different ESS technologies listed in section II.

## IV. TEST SYSTEMS

Proper test cases are key to show and implement the ideas of the researcher to interested readers. Some researchers and engineers with access to necessary power utility data like network and multi-period load data, are able to draw more interested readers since they have real-world systems as seen in approaches by [22], [29], [30] and [31]. If multi-period load data is not available, authors perform a snapshot power flow of a literature-based test system [46], [47], [48], [50] and [51] using the associated load data and perform their studies one period at a time or they scale the associated load data according to guidelines in the IEEE RTS-1996 [49] to develop a load curve that might be representative of the test system. These approaches were found in [12], [20], [33], [39], [40] and [44]. Reference [49] is the more frequently used literature-based test case as seen in approaches by [23], [24], [25], [41] and [43].

Some other authors in conjunction with power utilities or system operators attempt to scale obtained multi-period load data from those organizations to suit their selected literature-based case studies as seen in approaches by [26], [27], [32], [34], [41] and [45]. These literature-based test cases could also be modified based on assumptions that enable integration of the ESS and perhaps RES.

## V. CONCLUSION

In power system optimization considering ESSs, most researchers formulate generic ESS models independent of technology to capture their charge/discharge dynamics. BESS and PHS accounts for most of the models seen in literature. To our knowledge, there are currently no descriptive models capturing the characteristics of SMES and SC, hence providing a potential for future research. A generic ESS model would suffice for now for both technologies as well as other ESS technologies like CAES, FES and HES.

REFERENCES

[1] V.A. Boicea, "Energy Storage Technologies: The Past and the Present", Proc. of the IEEE, vol. 102, no. 11, pp. 1777-1794, Nov. 2014.
[2] X. Luo, J. Wang, M. Dooner and J. Clarke, "Overview of Current Development in Electrical Energy Storage Technologies and the Application Potential in Power System Operation", *App. Ener.*, 2015; 137:511-536.
[3] I.P. Gyuk and S. Eckroad, "EPRI-DOE Handbook of Energy Storage for Transmission and Distribution Applications", Electric Power Research Institute, Palo Alto, CA. 1001834, Dec. 2003.
[4] S. Rehman, L.M. Al-Hadhrami and Md.M. Alam, "Pumped Hydro Energy Storage System: A Technological Review", *Ren. & Sust. Ener. Rev.*, 2015; 44:586-598.
[5] G. Venkataramani, P. Parankusam, V. Ramalingam and J. Wang, "A Review on Compressed Air Energy Storage–A Pathway for Smart Grid and Polygeneration", *Ren. & Sust. Ener. Rev.*, 2016; 62:895-907.
[6] A.A.K. Arani, H. Karami, G.B. Gharehpetian and M.S.A. Hejazi, "Review of Flywheel Energy Storage Systems Structures and Applications in Power Systems and Microgrids", *Ren. & Sust. Ener. Rev.*, 2017; 69:9-18.
[7] K.C. Divya and J. Ostergaard, "Battery Energy Storage Technology for Power Systems–An Overview", *Elect. Pow. Syst. Res.*, 2009; 79:511-520

[8] W. Jing, C.H. Lai, S.H.W. Wong and M.L.D. Wong, "Battery-Supercapacitor Hybrid Energy Storage System in Standalone DC Microgrids: A Review", *IET Ren. Pow. Gen.*, vol. 11, iss. 4, pp. 461-469, May. 2017.

[9] A. Benato and A. Stoppato, "Pumped Thermal Electricity Storage: A Technology Overview", *Ther. Sci. & Eng. Prog.*, 2018; 6:301-315.

[10] X.Y. Chen, J.X. Jin, Y. Xin, B. Shu, CL. Tang et. al, "Integrated SMES Technology for Modern Power System and Future Smart Grid", *IEEE Trans. on App. Super.*, vol. 24, no. 5, pp. 3801605, Oct. 2014.

[11] G.L. Soloveichik, "Regenerative Fuel Cells for Energy Storage", Proc. of the IEEE, vol. 102, no. 6, pp. 964-975, Jun. 2014.

[12] A.S.A. Awad, D. Fuller, T.H.M. EL-Fouly and M.M.A. Salama, "Impact of Energy Storage Systems on Electricity Market Equilibrium", *IEEE Trans. on Sust. Ener.*, vol. 5, no. 3, pp. 875-885, Jul. 2014.

[13] R. Billinton and B. Bagen, "Reliability Considerations in the Utilization of Wind Energy, Solar Energy and Energy Storage in Electric Power Systems", *IEEE Proc. Int. Conf. on Prob. Meth. App. to Pow. Syst. – Stockholm, Sweden*, pp. 1-6, Jun. 2006.

[14] A. Ortega and F. Milano, "Generalized Model of VSC-Based Energy Storage Systems for Transient Stability Analysis", *IEEE Trans. on Pow. Syst.*, vol. 31, no. 5, pp. 3369-3380, Sep. 2016.

[15] M. Barnes, J. Kondoh, H. Asano, J. Oyarzabal, G. Ventakaramanan, R. Lasseter et. al, "Real-World Microgrids – An Overview", *IEEE Proc. Intl. Conf on Syst. of Syst. Eng. – San Antonio, USA*, pp. 1-8, Apr. 2007.

[16] H. Abdi, S.D. Beigvand and M. La Scala, "A Review of Optimal Power Flow Studies Applied to Smart Grids and Microgrids", *Ren. & Sust. Ener. Rev.*, 2017; 71:742-766.

[17] S.Y. Abujarad, M.W. Mustafa and J.J. Jamian, "Recent Approaches of Unit Commitment in the Presence of Intermittent Renewable Energy Sources: A Review", *Ren. & Sust. Ener. Rev.*, 2017; 70:215-223.

[18] D.K. Molzahn, F. Dorfler, H. Sandberg, S.H. Low, S. Chakrabarti, R. Baldick et. al, "A Survey of Distributed Optimization and Control Algorithms for Electric Power Systems", *IEEE Trans. on Smart Grid*, vol. 8, no. 6, pp. 2941-2962, Nov. 2017.

[19] M.B. Cain, R.P. O'Neill and A. Castillo, "History of Optimal Power Flow and Formulations", FERC Staff Technical Paper, Dec. 2012.

[20] S.F. Santos, D.Z. Fitiwi, M. Shafie-Khah, A. W. Bizuayehu, C.M.P. Cabrita et. al, "New Multistage and Stochastic Mathematical Model for Maximizing RES Hosting Capacity–Part I: Problem Formulation", *IEEE Trans. on Sust. Ener.*, vol. 8, no. 1, pp. 304-319, Jan. 2017.

[21] R.A. Jabr, S. Karaki and J.A. Korbane, "Robust Multi-Period OPF with Storage and Renewables", *IEEE Trans. on Pow. Syst.*, vol. 30, no. 5, pp. 2790-2799, Sept. 2015.

[22] P.S. Sauter, B.V. Solanki, C.A. Canizares, K. Bhattacharya and S. Hohmann, "Electric Thermal Storage System Impact on Northern Communities' Microgrids", *IEEE Trans. on Smart Grid*, [Early Access], Accepted Sept. 2017.

[23] M. Parvania, M. Fotuhi-Firuzabad and M. Shahidehpour, "Comparative Hourly Scheduling of Centralized and Distributed Storage in Day-Ahead Markets", *IEEE Trans. on Sust. Ener.*, vol. 5, no. 3, pp. 729-737, Jul. 2014.

[24] H. Pandzic, Y. Wang, T. Qiu, Y. Dvorkin and D.S. Kirschen, "Near-Optimal Method for Siting and Sizing of Distribution Storage in a Transmission Network", *IEEE Trans. on Pow. Syst.*, vol. 30, no. 5, pp. 2288-2300, Sep. 2015.

[25] N. Li and K.W. Hedman, "Enhanced Pumped Hydro Storage Utilization Using Policy Functions", *IEEE Trans. on Pow. Syst.*, vol. 32, no. 2, pp. 1089-1102, Mar. 2017.

[26] S. Bose, D.F. Gayme, U. Topcu and K.M. Chandy, "Optimal Placement of Energy Storage in the Grid", *51st IEEE Conf. on Dec. & Cont. – Hawaii, USA*, pp. 1-8, Dec. 2012.

[27] E. Sjodin, D.F. Gayme and U. Topcu, "Risk-Mitigated Optimal Power Flow for Wind Powered Grids", *IEEE Proc. Amer. Cont. Conf. – Montreal, Canada*, pp. 1-7, Jun. 2012.

[28] B.V. Solanki, A. Raghurajan, K. Bhattacharya and C.A. Canizares, "Including Smart Loads for Optimal Demand Response in Integrated Energy Management Systems for Isolated Microgrids", *IEEE Trans. on Smart Grid*, vol. 8, no. 4, pp. 1739-1748, Jul. 2017.

[29] Q. Li and V. Vittal, "Non-Iterative Enhanced SDP Relaxations for Optimal Scheduling of Distributed Energy Storage in Distribution Systems", *IEEE Trans. on Pow. Syst.*, vol. 32, no. 3, pp. 1721-1732, May. 2017.

[30] Q. Li and V. Vittal, "Convex Hull of the Quadratic Branch AC Power Flow Equations and Its Application in Radial Distribution Networks", *IEEE Trans. on Pow. Syst.*, vol. 33, no. 1, pp. 839-850, Jan. 2018.

[31] M. Nick, R. Cherkaoui and M. Paolone, "Optimal Siting and Sizing of Distributed Energy Storage Systems via Alternating Direction Method of Multipliers", *Elect. Pow. & Ener. Syst.*, 2015; 72:33-39.

[32] M. Nick, R. Cherkaoui and M. Paolone, "Optimal Allocation of Dispersed Energy Storage Systems in Active Distribution Networks for Energy Balance and Grid Support", *IEEE Trans. on Pow. Syst.*, vol. 29, no. 5, pp. 2300-2310, Sep. 2014.

[33] A. Gabash and P. Li, "Active-Reactive Optimal Power Flow in Distribution Networks with Embedded Generation and Battery Storage", *IEEE Trans. on Pow. Syst.*, vol. 27, no. 4, pp. 2026-2035, Nov. 2012.

[34] Q. Li, S. Yu, A.S. Al-Sumaiti and K. Turitsyn, "Micro Water-Energy Nexus: Optimal Demand-Side Management and Quasi-Convex Hull Relaxation", *arXiv preprint, arXiv: 1805.07626*, May 2018.

[35] M. Zidar, P.S. Georgilakis, N.D. Hatziargyriou, T. Capuder and D. Skrlec, "Review of Energy Storage Allocation in Power Distribution Networks: Applications, Methods and Future Research", *IET Gen., Trans. & Distr.*, vol. 10, no. 3, pp. 645-652, Mar. 2016.

[36] H. Saboori, R. Hemmati, S.M.S. Ghiasi and S. Dehghan, "Energy Storage Planning in Electric Power Distribution Networks – A State-of-the-Art Review", *Ren. & Sust. Ener. Rev.*, 2017; 79:1108-1121.

[37] Z. Wang, J. Zhong, D. Chen, Y. Lu and K. Men, "A Multi-Period Optimal Power Flow Model Including Battery Energy Storage", *IEEE PES General Meeting*, pp. 1-5, Jul. 2013.

[38] H. Daneshi and A.K. Srivastava, "Security-Constrained Unit Commitment with Wind Generation and Compressed Air Energy Storage", *IET Gen., Trans. & Distr.*, vol. 6, no. 2, pp. 167-175, Jan. 2012.

[39] Y.M. Atwa and E.F. El-Saadany, "Optimal Allocation of ESS in Distribution Systems with a High Penetration of Wind Energy", *IEEE Trans. on Pow. Syst.*, vol. 25, no. 4, pp. 1815-1822, Nov. 2010.

[40] D.E. Olivares, C.A. Canizares and M. Kazerani, "A Centralized Energy Management System for Isolated Microgrids", *IEEE Trans. on Smart Grid*, vol. 5, no. 4, pp. 1864-1875, Jul. 2014.

[41] D. Pozo, J. Contreras and E.E. Sauma, "Unit Commitment with Ideal and Generic Energy Storage Units", *IEEE Trans. on Pow. Syst.*, vol. 29, no. 6, pp. 2974-2984, Nov. 2014.

[42] R. Jiang, J. Wang and Y. Guan, "Robust Unit Commitment with Wind Power and Pumped Storage Hydro", *IEEE Trans. on Pow. Syst.*, vol. 27, no. 2, pp. 800-811, May. 2012.

[43] N. Li and K.W. Hedman, "Enhanced Pumped Hydro Storage Utilization Using Policy Functions", *IEEE Trans. on Pow. Syst.*, vol. 32, no. 2, pp. 1089-1102, Mar. 2017.

[44] A.A. Ibrahim, B. Kazemtabrizi, C. Bordin, C.J. Dent, J.D. McTigue and A.J. White, "Pumped Thermal Electricity Storage for Active Distribution Network Applications", *IEEE Proc. Manchester PowerTech Conf.*, pp. 1-6, Jun. 2017.

[45] S. Wogrin and D.F. Gayme, "Optimizing Storage Siting, Sizing and Technology Portfolios in Transmission-Constrained Networks", *IEEE Trans. on Pow. Syst.*, vol. 30, no. 6, pp. 3304-3313, Nov. 2015.

[46] Texas A&M University, "Electric Grid Test Case Repository", Dept. Electrical & Computer Engineering, 2016. [Online]. Available: https://electricgrids.engr.tamu.edu/electric-grid-test-cases/.

[47] M.E. Baran and F.F. Wu, "Network Reconfiguration in Distribution Systems for Loss Reduction and Load Balancing", *IEEE Trans. on Pow. Del.*, vol. 4, no. 2, pp. 1401-1407, Apr. 1989.

[48] K. Rudion, A. Orths, Z.A. Styczynski and K. Strunz, "Design of Benchmark of Medium Voltage Distribution Network for Investigation of DG Integration", *IEEE PES General Meeting*, pp. 1-6, Oct. 2006.

[49] C. Grigg, P. Wong, P. Albrecht, R. Allan, M. Bhavaraju, R. Billinton et. al, "The IEEE Reliability Test System–1996", *IEEE Trans. on Pow. Syst.*, vol. 14, no. 3, pp. 1010-1020, Aug. 1999.

[50] M.E. Baran and F.F. Wu, "Optimal Capacitor Placement on Radial Distribution Systems", *IEEE Trans. on Pow. Del.*, vol. 4, no. 1, pp. 725-734, Jan. 1989.

[51] IEEE PES AMPS DSAS Test Feeder Working Group, "Test Feeder Cases". [Online]. Available: http://sites.ieee.org/pes-testfeeders/resources/.